\documentclass[11pt]{article}
\usepackage{graphicx}
\usepackage{amsfonts}
\usepackage{theorem}
\newtheorem{Lem}{Lemma}[section]

\newtheorem{The}[Lem]{Theorem}

\newtheorem{Cor}[Lem]{Corollary}

\newtheorem{Rem}[Lem]{Remark}

\newcommand{\qed}{\hbox{\rule{6pt}{6pt}}}

\begin{document}
\title{Refinements of the trace inequality of Belmega, Lasaulce and Debbah}
\author{Shigeru Furuichi$^1$\footnote{E-mail:furuichi@chs.nihon-u.ac.jp} and Minghua Lin$^2$\footnote{E-mail:lin243@uregina.ca}\\
$^1${\small Department of Computer Science and System Analysis,}\\
{\small College of Humanities and Sciences, Nihon University,}\\
{\small 3-25-40, Sakurajyousui, Setagaya-ku, Tokyo, 156-8550, Japan}\\
$^2${\small Department of Mathematics and Statistics,}\\
{\small University of Regina, Regina, Saskatchewan, Canada S4S 0A2}}
\date{}
\maketitle

{\bf Abstract.} In this short paper, we show a certain matrix trace inequality and then
give a refinement of the trace inequality proven by Belmega, Lasaulce and Debbah.
In addition, we give an another improvement of their trace inequality.
\vspace{3mm}


\vspace{3mm}
\vspace{3mm}


\section{Introduction}
Recently, E.-V.Belmega, S.Lasaulce and M.Debbah obtained the following elegant trace inequality for positive definite matrices.

\begin{The} \label{the01} {\bf (\cite{BSL})}
 For  positive definite matrices $A, B$ and positive semidefinite matrices $C,D$, we have
\begin{equation}  \label{ineq01}
Tr[(A-B)(B^{-1}-A^{-1})+(C-D)\left\{(B+D)^{-1}-(A+C)^{-1}\right\}] \geq 0.
\end{equation}
\end{The}
In this short paper, we first prove a certain trace inequality for products of matrices, and then as its application, we give a simple proof of (\ref{ineq01}).
At the same time, our alternative proof gives a refinement and of Theorem \ref{the01}. An another improvement of the Theorem \ref{the01} is also considered at the end of the paper.

\section{Main results}
In this section, we prove the following theorem.

\begin{The} \label{the02}
For  positive definite matrices $A, B$ and positive semidefinite matrices $C,D$, we have
\begin{eqnarray}
&&Tr[(A-B)(B^{-1}-A^{-1})+(C-D)\left\{(B+D)^{-1}-(A+C)^{-1}\right\}] \nonumber \\
&&\geq \vert Tr[(C-D)(B+D)^{-1}(A-B)(A+C)^{-1}]\vert . \label{ineq02}
\end{eqnarray}
\end{The}

To prove this theorem, we need a few lemmas.

\begin{Lem} {\bf (\cite{BSL})}  \label{lem01}
For positive definite matrices $A, B$ and positive semidefinite matrices $C,D$, and Hermitian matrix $X$, we have
$$
Tr[XA^{-1}XB^{-1}] \geq Tr[X(A+C)^{-1}X(B+D)^{-1}].
$$
\end{Lem}


\begin{Lem} \label{lem02}
For any matrices $X$ and $Y$, we have
$$
Tr[X^*X]+ Tr[Y^*Y] \geq 2\vert Tr[X^*Y]\vert.
$$
\end{Lem}
{\it Proof}:
Since $Tr[X^*X] \geq 0$, by the fact that the arithmetical mean is greater than the geometrical mean and  Cauchy-Schwarz inequality, we have
$$
\frac{Tr[X^*X]+Tr[Y^*Y]}{2} \geq \sqrt{Tr[X^*X] Tr[Y^*Y] } \geq  \vert Tr[X^*Y]\vert .
$$

\hfill \qed

\begin{The} \label{the03}
For Hermitian matrices  $X_1, X_2$ and positive semidefinite matrices $S_1, S_2$, we have
$$
Tr[X_1S_1X_1S_2]+Tr[X_2S_1X_2S_2] \geq 2 \vert Tr[X_1S_1X_2S_2]\vert.
$$
\end{The}
{\it Proof}:
Applying Lemma \ref{lem02}, we have
\begin{eqnarray*}
&&Tr[X_1S_1X_1S_2]+Tr[X_2S_1X_2S_2] \\
&&\hspace{5mm} = Tr[(S_2^{1/2}X_1S_1^{1/2})(S_1^{1/2}X_1S_2^{1/2})]+Tr[(S_2^{1/2}X_2S_1^{1/2})   (S_1^{1/2}X_2S_2^{1/2})]\\
&&\hspace{5mm} \geq 2 \vert Tr[(S_2^{1/2}X_1S_1^{1/2})(S_1^{1/2}X_2S_2^{1/2}) ]\vert \\
&&\hspace{5mm} = 2 \vert Tr[X_1S_1X_2S_2]\vert.
\end{eqnarray*}
\hfill \qed

\begin{Rem}
Theorem \ref{the03} can be regarded as a kind of the generalization of Proposition 1.1 in \cite{FKY}.
\end{Rem}

{\it Proof of Theorem \ref{the02}}:
By Lemma \ref{lem01}, we have
\begin{eqnarray*}
Tr[(A-B)(B^{-1}-A^{-1})]&=&Tr[(A-B)B^{-1}(A-B)A^{-1}]\\
&\geq & Tr[(A-B)(A+C)^{-1}   (A-B)(B+D)^{-1}]\\
&=& Tr[(A-B)(B+D)^{-1}(A-B)(A+C)^{-1} ].
\end{eqnarray*}
Thus the left hand side of the inequality (\ref{ineq02}) can be bounded from below:
\begin{eqnarray}
&&Tr[(A-B)(B^{-1}-A^{-1})+(C-D)\left\{(B+D)^{-1}-(A+C)^{-1}\right\}] \nonumber \\
&& \geq Tr[(A-B)(B+D)^{-1}(A-B)(A+C)^{-1}+(C-D)(B+D)^{-1}(C-D)(A+C)^{-1}] \nonumber \\
&& \hspace{5mm} +Tr[(C-D)(B+D)^{-1}(A-B)(A+C)^{-1}] \nonumber \\
&& \geq 2 \vert Tr[(C-D)(B+D)^{-1}(A-B)(A+C)^{-1}]  \vert\nonumber \\
&&  \hspace{5mm} +Tr[(C-D)(B+D)^{-1}(A-B)(A+C)^{-1}]  \label{ineq03}
\end{eqnarray}
Throughout the process of the above, Theorem \ref{the03} was used in the second inequality.
Since we have the following equation,
\begin{eqnarray*}
&&Tr[(C-D)(B+D)^{-1}(A-B)(A+C)^{-1}] \\
&& \hspace{5mm} =Tr[(C-D)(B+D)^{-1}]- Tr[(C-D)(A+C)^{-1}]\\
&&  \hspace{7mm} -Tr[(C-D)(B+D)^{-1}(C-D)(A+C)^{-1}]
\end{eqnarray*}
we have $ Tr[(C-D)(B+D)^{-1}(A-B)(A+C)^{-1}] \in \mathbb{R}$.
Therefore  we have
$$
(\ref{ineq03}) \geq \vert Tr[(C-D)(B+D)^{-1}(A-B)(A+C)^{-1}]\vert .
$$

\hfill \qed

\section{An another improvement of the inequality (\ref{ineq01})}
In this section, we show the following trace inequality.
\begin{The}  \label{the31}
For  positive definite matrices $A, B$ and positive semidefinite matrices $C,D$, we have
\begin{equation}   \label{ineq31}
Tr[(A-B)(B^{-1}-A^{-1})+4(C-D)\left\{(B+D)^{-1}-(A+C)^{-1}\right\}]\geq 0.
\end{equation}
\end{The}

To prove this theorem, we use the following lemmas, which are proven by the similar way of Lemma \ref{lem02} and Theorem \ref{the03} in the previous section.
\begin{Lem} \label{lem32}
For any matrices $X$ and $Y$, any positive real numbers $a$ and $b$, we have
$$
a\cdot Tr[X^*X]+ b \cdot Tr[Y^*Y] \geq 2 \sqrt{ab}\cdot  \vert Tr[X^*Y]\vert.
$$
\end{Lem}
Applying this lemma, we have the following lemma.
\begin{Lem} \label{lem33}
For Hermitian matrices  $X_1, X_2$, positive semidefinite matrices $S_1, S_2$ and any positive real numbers  $a$ and $b$, we have
$$
a\cdot Tr[X_1S_1X_1S_2]+ b \cdot Tr[X_2S_1X_2S_2] \geq  2 \sqrt{ab} \cdot \vert Tr[X_1S_1X_2S_2]\vert.
$$
\end{Lem}

{\it Proof of Theorem \ref{the31}}:
By the similar way to the proof of Theorem \ref{the02}, applying Lemma \ref{lem32} as $a=1$ and $b=4$,
the left hand side of the inequality of (\ref{ineq31}) can be bounded from the below:
\begin{eqnarray}
&&Tr[(A-B)(B^{-1}-A^{-1})+4(C-D)\left\{(B+D)^{-1}-(A+C)^{-1}\right\}] \nonumber \\
&& \geq Tr[(A-B)(B+D)^{-1}(A-B)(A+C)^{-1}+4(C-D)(B+D)^{-1}(C-D)(A+C)^{-1}] \nonumber \\
&& \hspace{5mm} +Tr[4(C-D)(B+D)^{-1}(A-B)(A+C)^{-1}] \nonumber \\
&& \geq 4 \vert Tr[(C-D)(B+D)^{-1}(A-B)(A+C)^{-1}]  \vert\nonumber \\
&&  \hspace{5mm} +4\cdot Tr[(C-D)(B+D)^{-1}(A-B)(A+C)^{-1}]  \geq 0, \nonumber
\end{eqnarray}
since  $ Tr[(C-D)(B+D)^{-1}(A-B)(A+C)^{-1}] \in \mathbb{R}$.
\hfill \qed

\begin{Rem}
Here we note that we have $Tr[(A-B)(B^{-1}-A^{-1})]\geq 0$.
However we have the possibility that $Tr[(C-D)\left\{(B+D)^{-1}-(A+C)^{-1}\right\}]$  takes a negative value.
Therefore Theorem \ref{the31} is an improvement of Theorem \ref{the01}.
\end{Rem}

\begin{Cor} \label{cor35}
For  positive definite matrices $A, B$, positive semidefinite matrices $C,D$ and positive real number $r$, we have
\begin{equation}   \label{ineq32}
Tr[(A-B)(B^{-1}-A^{-1})+4(C-D)\left\{(r B+D)^{-1}-(r A+C)^{-1}\right\}]\geq 0.
\end{equation}
\end{Cor}
{\it Proof}:
Put $A = r A_1$ and $B=r B_1$ for  positive definite matrices $A_1$ and $B_1$,  in Theorem \ref{the31}.
\hfill \qed

\begin{Rem}
In the case of $r=2$ in Corollary \ref{cor35}, the inequality (\ref{ineq32}) corresponds to the scalar inequality:
$$
(\alpha-\beta)\left(\frac{1}{4\beta}-\frac{1}{4\alpha}\right)+\left(\gamma-\delta\right)\left(\frac{1}{2\beta+\delta}-\frac{1}{2\alpha+\gamma}\right)\geq 0
$$
for positive real numbers $\alpha$ and $\beta$, nonnegative real numbers $\gamma$ and $\delta$.
\end{Rem}

\end{document}